\documentclass[11pt]{amsart}

\usepackage[utf8]{inputenc}
\usepackage[T1]{fontenc}

\usepackage[english]{babel}
\usepackage{fullpage}
\usepackage{enumitem}

\usepackage{amssymb,libertine}
\usepackage{amssymb}
\usepackage{array}
\usepackage{bbm}

\usepackage{mathtools}
\usepackage[svgnames]{xcolor}
\usepackage{hyperref}
\usepackage{comment}
\usepackage[colorinlistoftodos]{todonotes}

\usepackage{mathrsfs}
\usepackage{wrapfig}
\usepackage{color}
\usepackage{tikz}

\usepackage{graphicx}
\usepackage{epstopdf}

 \usepackage{latexsym}
 \usepackage{amsfonts} 

\usepackage{amsmath}
\usepackage{amsthm}
\usepackage[export]{adjustbox}
\def\lap{\mathscr{L}}
   \def\CC{\mathbb{C}}
    \def\DD{\mathbb{D}}
    \def\NN{\mathbb{N}}
    
    \def\RR{\mathbb{R}}
    
    \newtheorem{Proposition}{Proposition}
\newtheorem{Theorem}[Proposition]{Theorem}

\newtheorem{Definition}[Proposition]{Definition}

\newtheorem{Remark}[Proposition]{Remark}
\newtheorem{Note}[Proposition]{Note}
  
\def\be{\begin{equation}}
\def\ee{\end{equation}}

\def\le{\leqslant}

\def\bd{\begin{Definition}}
\def\ed{\end{Definition}}
\def\bt{\begin{Theorem}}
\def\et{\end{Theorem}}

 \def\NDP{ [n_j,m_j]_{j\in \NN}}

\def\Om{\Omega}

\def\epsilon{\varepsilon}
\def\bel{\begin{equation}\label}
\def\ee{\end{equation}}

\def\Ei\text{Ei}

\def\phi{\varphi}

\usepackage{graphicx} 
\graphicspath{{out/}{./}}

\title{Approximations of Functions With Essential Singularities with Applications to Painlevé's First Transcendent} 
\author{Nicholas Castillo  }
\date{July 2023}

\begin{document}

\maketitle

\section{Abstract}
In this work we develop an algorithmic procedure for associating a function defined on the Riemann surface of the $\log$ to given asymptotic data from a function at an essential singularity.  We do this by means of rational approximations (Padé approximants) used in tandem with Borel-Écalle summation.  Our method is capable of handling situations where classical methods either do not work or converge very slowly eg. \cite{DunLutz}.  We provide a general outline of the procedure and then apply it to generating approximate tritronquée solutions to Painlevé's first equation ($\text{P}_\text{I}$).  Our approximations (including  $\text{P}_\text{I}$) are written as a finite linear combination of exponential integrals $\text{Ei}^+$. Furthermore, from \cite{GRA} (the relevant results are summarized in \S\ref{DyEiSec}) we have explicit rational approximations for each $\text{Ei}^+$ and thus for the approximation as a whole. In addition to rational approximations of $\text{P}_\text{I}$, we provide the first hundred or so poles of a tritronquée solution with essentially arbitrary  accuracy which is dependent upon the order of Padé used. 

\section{The General Method}\label{GenMet}
We begin with a situation frequently encountered in the practice of solving ODE, PDE or physical problems of various flavors.  Given a finite set of perturbative data in the form of an asymptotic power series

\begin{equation}\label{APS}
    \sum_{k=0}^{2N}\frac{a_k}{x^{k+1}}
\end{equation}
we will construct an approximate solution to the problem at hand 
 having the first $2N$ terms of its asymptotic power series given by \eqref{APS}. 

\subsection{Our Procedure}

\begin{enumerate}\label{Procedure}
    \item Apply the Borel transform $\mathscr{B}$ to \eqref{APS} and obtain the polynomial: 
\begin{equation}\label{Borel}
   F_{2N}(p)= \mathscr{B}\left(\sum_{k=0}^{2N}\frac{a_k}{x^{k+1}}\right)=\sum_{k=0}^{2N}\frac{a_k p^k}{k!}
\end{equation}
    \item Develop a Padé approximant $\mathcal{P}_m^n(p)$ of \eqref{Borel} about the origin; (in our example we use diagonal Padé so $n=m=N$)\\
    \item Break  $\mathcal{P}_m^n(p)$ into partial fractions with simple poles located at the points $p_k$.\\
    \item Continuing to follow the procedure of Borel summation, we take the Laplace transform $\lap_\theta$ along a direction $\theta$ disjoint from the set of poles $\{p_k\}_{k=1}^N$.  
\end{enumerate}

The directional Laplace transform $\lap_\theta$ is a bounded operator $L_\nu^1\left ({\rm{e}}^{i\theta}\RR^+\right)$\footnote{$L_\nu^1\left ({\rm{e}}^{i\theta}\RR^+\right)$ is a space of exponentially weighted $L^1$ functions; $\{F:{\rm{e}}^{i\theta}\RR^+\to \CC \big \vert F(p){\rm{e}}^{-\nu p}\in L^1\left ({\rm{e}}^{i\theta}\RR^+\right)\}$ }$\to\mathscr{O}\left(\{\Re (x{\rm{e}}^{i\theta})>\nu\}\right)$ defined by the contour integral along the ray $\{\arg p=\theta\}$
\begin{equation}
    \lap_\theta [F](x)=\int_0^{\infty {\rm{e}}^{i\theta}}F(p){\rm{e}}^{-px}dp
\end{equation}
We are left with the function 
\begin{equation}
    f(x)=\lap_\theta [F_{2N}](x)= \sum_{k=0}^{N} c_k\lap_\theta \left\{\frac{1}{p-p_k}\right\}(x)
\end{equation}
with domain of analyticity $\Pi_\theta \coloneqq\{x\in\CC:\Re (x{\rm{e}}^{i\theta})>0\}$.  Each pole $p_k$ corresponds to a Stokes direction, we let $\theta_k=\arg p_k$ 
and pick $\theta_1,\theta_2$ to be the two neighboring directions of the contour of integration $\{t{\rm{e}}^{i\theta}\}_{t\geq0}$; ie. $\theta_1<\theta<\theta_2$.  Each $\lap_\theta \left\{\frac{1}{p-p_k}\right\}$ may be written in terms of the representation of the exponential integral $\mathrm{Ei}^+$ we have chosen to work with defined by \eqref{LiEiRatAp}.
\begin{equation}
    \lap_\theta \left\{\frac{1}{p-p_k}\right\}=-{\rm{e}}^{-p_k x}\mathrm{Ei}^+\left(p_k x\right)
\end{equation}
Therefore we have the following representation 

\begin{equation}\label{EiAppr}
    f_{N,\theta}(x)=\lap_\theta [F_{2N}](x)= -\sum_{k=0}^{N} c_k {\rm{e}}^{-p_k x}\mathrm{Ei}^+\left(p_k x\right)
\end{equation}

Furthermore, each $\text{Ei}^+$ admits a factorial expansion given by \eqref{eq:EidsumP1}. We may arrange for all cuts defining the boundary of convergence for each expansion to coincide by choosing an appropriate parameter $\beta$; c.f. \cite{GRA} .

We may then analytically continue the function element $(f_{N,\theta},\Pi_\theta)$ to the domain 
\begin{equation}
    \mathscr{D}(\theta_1,\theta_2)\coloneqq \bigcup_{\varphi\in (\theta_1,\theta_2)} \Pi_\varphi
\end{equation}
a sectorial neighborhood of infinity.  \\

 \begin{Remark}
 The problem in which one knows the nature of the Riemann surface associated to a given Maclaurin series has been studied and laid to rest in \cite{Unif}.  In this case, the uniformization map and its inverse for the Riemann surface used together with the truncated Maclaurin series (as the given data in place of our asymptotic power series \eqref{APS})  yields the best possible approximation to the actual function.
 \end{Remark}

\section{Painlevé's First Equation}
\subsection{Tritronqué solutions }\label{triSol}
In the present work, we use the normalization in which $\text{P}_\text{I}$ takes the form 
\begin{equation}
    \label{P1}
    y''=6y^2-z
\end{equation}

\begin{figure}\label{LogErrRHP}
  \centering 
\includegraphics[scale=0.4]{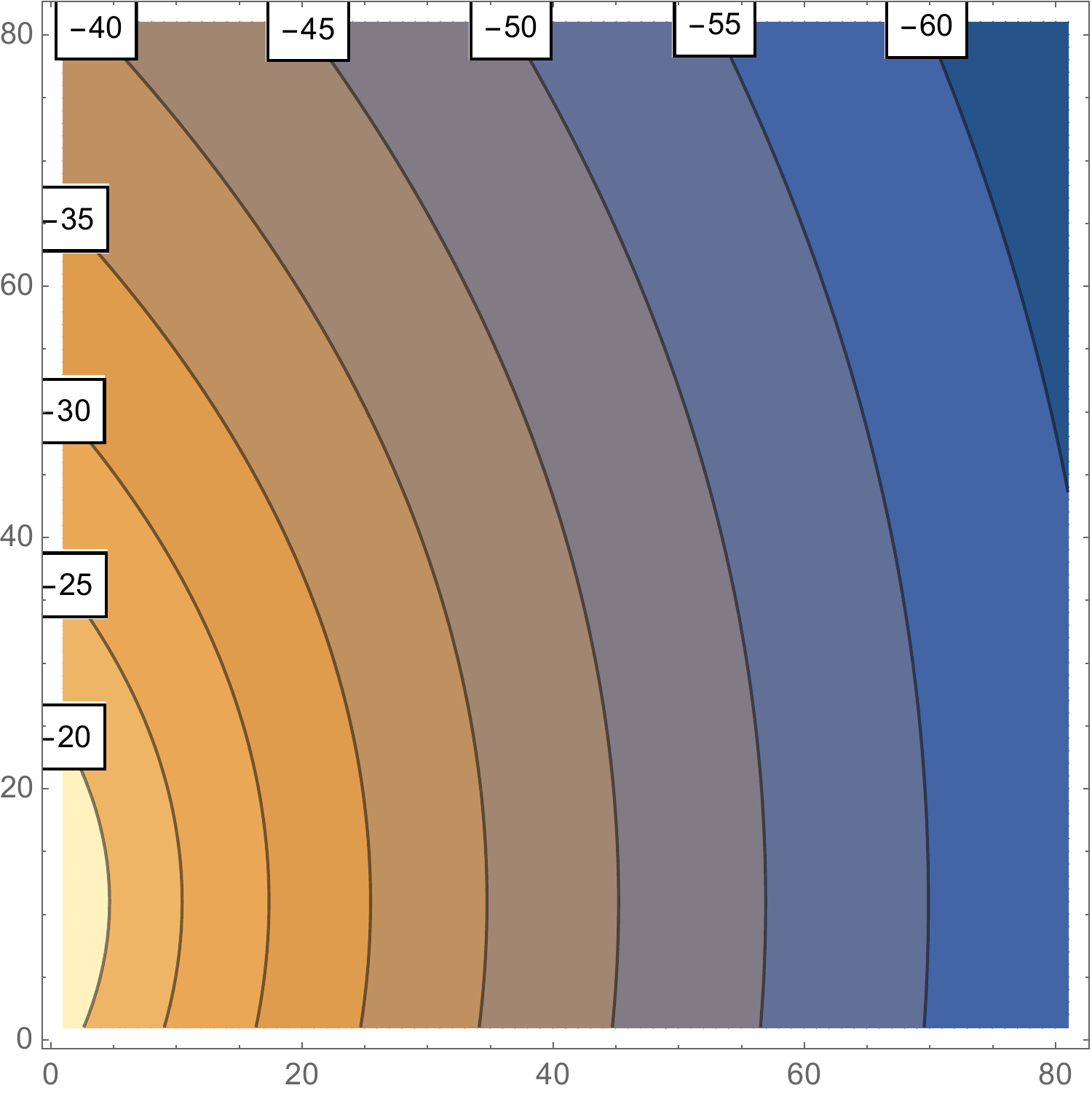} 
\caption{Modulus of error log-plot generated from a fifty exponential integral approximation.  The $x,y$ axes are the number of steps of size $\frac{1}{10}$ in the real and imaginary directions respectively starting at $1-i$. }
\end{figure}

\begin{figure}\label{LowErrLHP}
  \centering 
\includegraphics[scale=0.4]{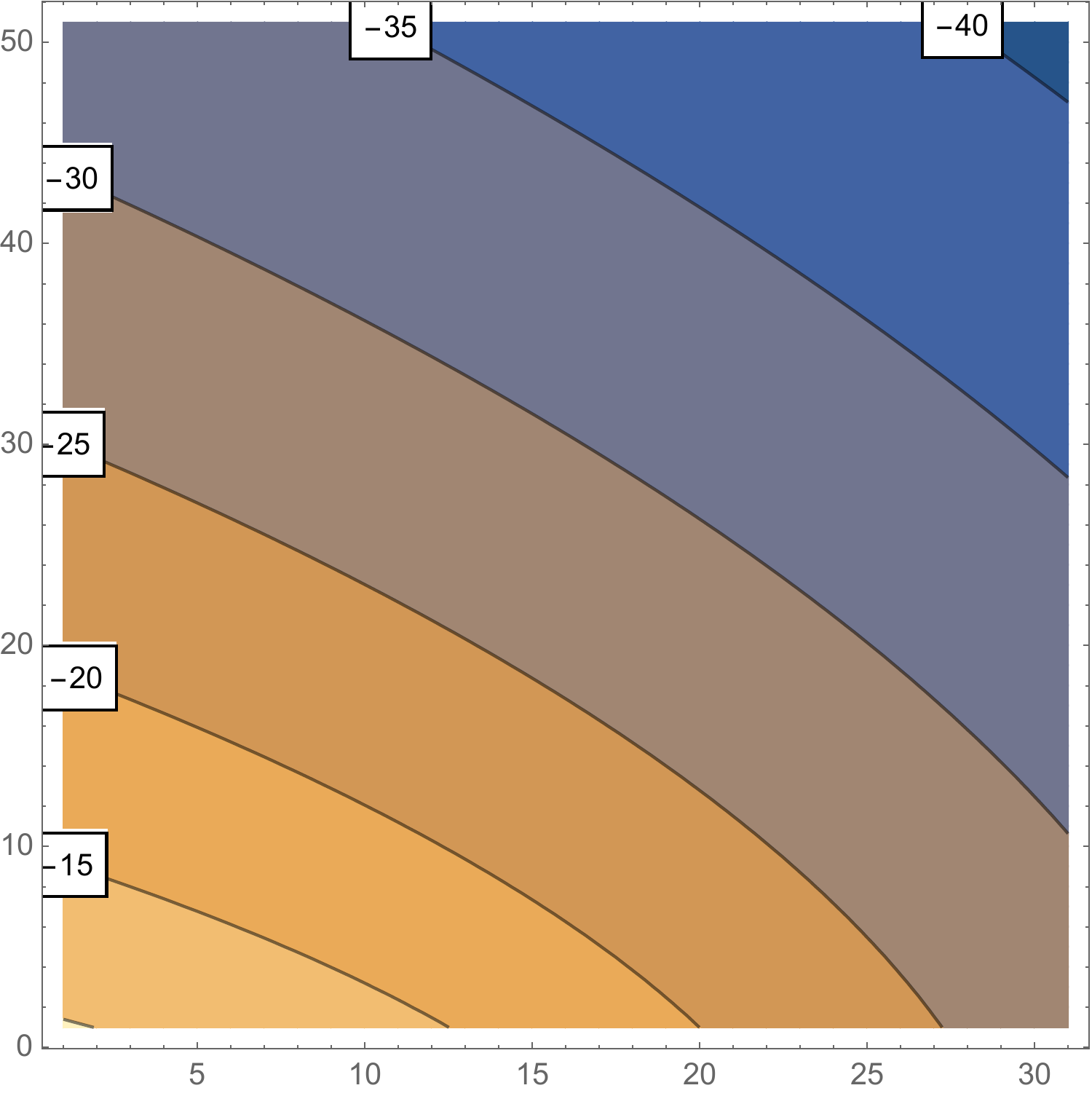} 
\caption{Modulus of error log-plot generated from a fifty exponential integral approximation.  The $x,y$ axes are the number of steps of size $\frac{1}{10}$ in the real and imaginary directions respectively starting at $-1-\frac{3}{2}i$. }
\end{figure}

From here, we use coordinates which were inspired by Boutroux's original work on $\text{P}_\text{I}$ \cite{Boutroux}

\begin{equation}
    \label{CoordChng}
    x=\frac{\left(24z\right)^{5/4}}{30} \quad , \quad y=-\sqrt{\frac{z}{6}}\left(1+h(x)\right)
\end{equation}

and \eqref{P1} now reads 

\begin{equation}
\label{hEqu}
    h''+\frac{h'}{x}+h-\frac{4}{25x^2}+\frac{1}{2}h^2-\frac{4}{25x^2}h=0
\end{equation}
A reader familiar with resurgence will recognize \eqref{CoordChng} as the Écalle critical time associated with \eqref{P1} which ensures the equation takes a form suitable for Borel-Écalle summation.  Taking the inverse Laplace transform we arrive at the following convolution equation defined on the Borel plane: 

\begin{equation}
    \left(p^2+1\right)H(p)-1*pH(p)-\frac{4}{25}p+\frac{1}{2}H*H-\frac{4}{25}p*H=0
\end{equation}
where $H(p)=\lap^{-1}h$.  The polynomial coefficient $p^2+1$ and non-linearity of \eqref{P1} tells us to expect two arrays of singular points $\pm i \NN$ in the Borel p-plane, each of which corresponds to a Stokes direction.  For the general theory, we refer the reader to \cite{IMRN} or for even more generality \cite{Duke}.  Following our procedure described in \S\ref{GenMet}, we first generate a truncated asymptotic power series about infinity from \eqref{hEqu}. Note that the full series is necessarily (factorially) divergent.  The first few terms of which are 

\begin{multline}\label{AsympExpP1}
   \frac{4}{25 x^2}-\frac{392}{625
   x^4}+\frac{6272}{625
   x^6}-\frac{141196832}{390625
   x^8}+\frac{9039055872}{390625
   x^{10}}\\-\frac{565008634278144}{24
   4140625
   x^{12}}+ \frac{81365672232484864}{244140625
   x^{14}}+\cdots
\end{multline}
   Continuing to follow \S\ref{GenMet} we truncate \eqref{AsympExpP1} to order $2N$ and compute the Borel transform $\mathscr{B}$ according to \eqref{Borel} from which we obtain a truncated Maclaurin series in the Borel $p$-plane.  Using the truncated series we develop a $[N,N]$ Padé approximant about $p=0$ and perform partial fraction decomposition.  It is this decomposition which produces the data set $\{c_k,p_k\}_{k=1}^N$ where $c_k$ are the residues corresponding to the poles $p_k$.  It is exactly this set which generates our tritronquée approximation according to \eqref{EiAppr} when we return to the physical domain by a Laplace transform along any direction $\theta\neq \pi/2 \mod \pi$.  We refer the reader to Figures \ref{LogErrRHP} and \ref{LowErrLHP} which show contours of the modulus of error log-plot generated from a fifty exponential integral approximation.  The errors were computed by applying the nonlinear differential operator defined by the left side of \eqref{hEqu} to the Ei approximation and evaluating it over a grid comprised of $10^3$ points.  
   
   \subsection{Tritronquée Pole Locations}
   We use the approximation obtained in \S\ref{triSol} to provide the two free coefficients $c_0=y(z_0), c_1=y'(z_0)$ found in a power series expansion of \eqref{P1} about a point $z_0$.  In our calculations we took $z_0=5+i$.  We develop a power  from relations computed through \eqref{P1}
 \begin{equation}
      c_0=y(z_0), \quad c_1=y'(z_0), \quad c_2=-\frac{1}{2}z_0+3c_0^2, \quad c_3=2c_1c_0-\frac{1}{6}
 \end{equation}
     with a general term: 
     \begin{equation}\label{GenTermSer}
            c_n=\frac{6}{n(n-1)}\sum_{j=0}^{n-2}c_j c_{n-2-j}, \quad T(z)=\sum_{n=0}^{2N} c_n(z-z_0)^n
     \end{equation}

\begin{figure}\label{fig:PolePlot}
  \centering 
\includegraphics[scale=0.4]{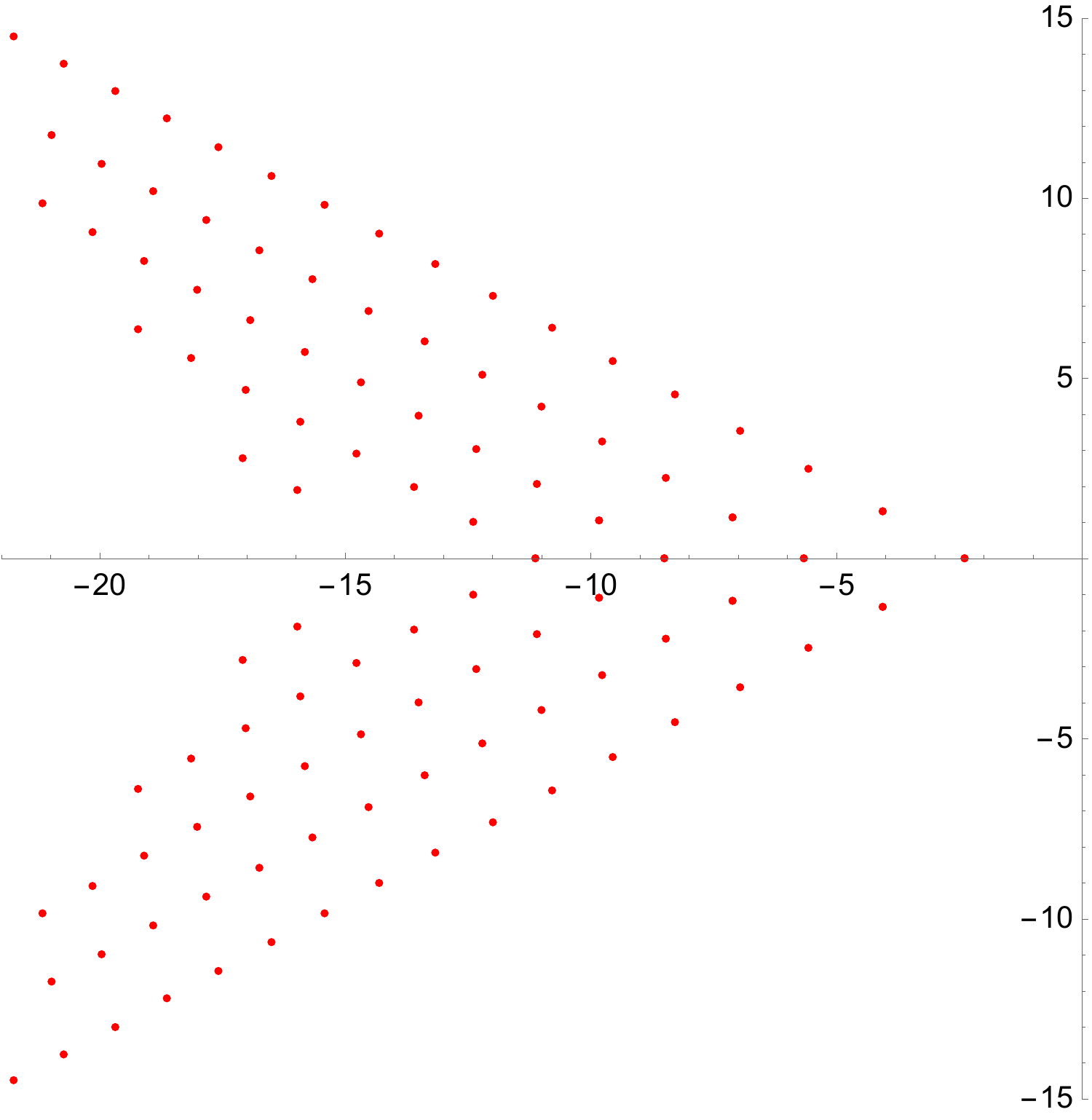} 
\caption{Pole locations of a tritronquée solution to $\text{P}_\text{I}$. }
\end{figure}

   From \eqref{GenTermSer} we compute a Padé approximant based at $z_0$.  It is this rational function which will provide our approximation for the pole locations of a tritronquée solution to $\text{P}_\text{I}$.   However, it is well known that Padé will often create \emph{spurious poles} which are unrelated to the function one is studying \cite{Stahl2}.  To handle this issue we compute the residue of each pole and distinguish between spurious and non-spurious by the magnitude of the residues.  The size of the residues of genuine poles tend to be significantly larger in magnitude than most spurious poles \cite{Gonnet}.    The exact transition from spurious to non-spurious is not precise and one should check borderline cases by using computer plotting software.  From this process we have plotted some of the poles for tritronquée in Figure \ref{fig:PolePlot}.

\section{Estimate of error}
Following the results from Stahl's foundational paper \cite{Stahl} we compute the error in our approximation.  We use the notation $\widehat{\CC}$ to denote the Riemann sphere.  As discussed in \S\ref{PadeSection}, Padé approximants are only guaranteed to converge in a potential theoretic sense i.e. convergence in capacity.  However, in our case the compliment of the minimal capacitor is the simply connected domain $\Om=\widehat{\CC}\setminus\{it:|t|\geq1\} $.  This makes the calculation of the associated greens function significantly more straightforward c.f. \cite{Stahl}.  In this case the Green's function is defined by

\begin{equation}
    g_\Om(z;0)=-\log |\psi(z)|
\end{equation}
where $\psi:\Om\to \DD$ is the conformal map guaranteed by the Riemann mapping theorem with normalization $\psi(0)=0$ defined by:

\begin{equation}
    \psi(w)=i\left(\frac{1-\sqrt{1+w^2}}{w}\right)
\end{equation}
and therefore 
\begin{equation}
    g_\Om(w;0)=-\log\left|\frac{1-\sqrt{1+w^2}}{w}\right|
\end{equation}
\newpage
Therefore, according to \eqref{CapG} the Green's function will be 
\begin{equation}
   G_\Omega(w)= \left|\frac{1-\sqrt{1+w^2}}{w}\right|
\end{equation}
and we have the potential theoretic estimate 
\bt
Let $h(x)$ be a solution to \eqref{hEqu}. 
Then for every $\varepsilon>0$ and compact $V \subset \widehat{\CC}\setminus\{it:|t|\geq1\}$
    $$\lim_{j\to \infty} cap\left\{x\in V: \left|h(x)-[n_j,m_j](x)\right|>\left(\left|\frac{1-\sqrt{1+x^2}}{x}\right|+\varepsilon\right)^{n_j+m_j}\right\}=0$$
   where  $\left \{[n_j,m_j]\right\}_{j\in\NN}$ is a near diagonal Padé approximant. 
\et
\begin{proof}
    This is immediate from Theorem 1.2(i) in Stahl's paper  \cite{Stahl}.   
\end{proof}

\section{Appendix}
\subsection{Padé Approximation}\label{PadeSection}

Approximation theory has a long and rich history touching many different fields of mathematics.  Most people first encounter such principles when studying Taylor series and families of orthogonal polynomials.  However, the use of polynomial approximation is not always the most efficient when dealing with problems of natural origin, especially when singularities are involved.  There can be issues with rates of convergence as well as divergence of the given expansion.  Additionally, sometimes it's only practical to extract a few terms of a given expansion which may not yield the desired accuracy.  In both the convergent and divergent cases we have the tools of Padé approximants.  Such tools enable one to probe for the presence of poles and branch points with surprising accuracy. 

\begin{Note}\label{PadAssumptionsNote}
The functions considered here will be holomorphic at infinity with a set $E$ of singularities of zero capacity. This is equivalent to $f$ having (possibly multivalued) analytic continuations along any path starting at $\infty$ and contained in $\overline{\CC}\setminus E$ \cite{Stahl}.   
\end{Note}

\bd[Padé approximant \cite{Stahl}]
For a given function $f$ analytic at infinity we define the $[n/m]$ Padé approximant as
the rational function $P_m^n(z^{-1})=A_n(z^{-1})/B_m(z^{-1})$ where $A_n,B_m$  are polynomials in $z^{-1}$ of degrees at most $n$ and $m$ respectively such that 
$$f(z^{-1})-P_m^n(z^{-1})=O(z^{-n-m-1})$$
as $z\to \infty$.  If we normalize so that the constant term in the denominator is equal to $1$ then such rational functions are unique for a given $f$ and $n,m\in \NN$.  The sequence $\{[n/n]\}_{n\in \NN}$ is called diagonal and $\{[n_j/m_j]\}_{j\in \NN}$ such that $m_j\to \infty$ and $n_j/m_j\to 1$ as $j\to \infty$ are called near-diagonal. 

\ed

We remark that the first $n+m+1$ terms of the Taylor expansion for $P_m^n(z)$ about infinity match the first $n+m+1$ terms of the expansion for $f$.  Additionally, we must note one drawback; convergence is not typically pointwise but rather in the sense capacity.  Such phenomena will be described below.   

These ideas go back to the work of Georg Frobenius and of course Henri Padé but work in this area has continued well into the twentieth and twenty first centuries. We begin with one important result attributed to Nuttall and Pommerenke:

\bt[Nuttall-Pommerenke theorem \cite{Nuttall}, \cite{Pom}]
Given a function $f$ satisfying the conditions of Note \ref{PadAssumptionsNote} such that $f$ is single valued on $\CC\setminus E$.  Then for any $\varepsilon>0, \lambda\in (0,1)$ and $V \subset \subset \CC\setminus E$ we have 
$$\lim_{n,m\to \infty} cap\{z\in V: |f(z)-P_m^n(z)|>\varepsilon^{n+m}\}=0$$
with $(n,m)\in \NN^2$ satisfying 
$$\lambda n\leq m\leq \frac{n}{\lambda}$$
for $n,m\to \infty$. 
\et
The above is an example of convergence in capacity which is analogous to that of convergence in measure.  We refer the reader  to \cite{SaffLog} for a full treatment of basic potential theory. 


What is remarkable is that the Padé sequence provides a maximal domain of single-valuedness and convergence in capacity (up to sets of capacity zero). Such a domain emerges so that the boundary of convergence is of minimal logarithmic capacity \cite{Stahl}, \cite{Stahl2}.  Once more, the location of the poles in the Padé sequence cluster around the boundary. We now state two result of Stahl  which will be relevant.  

\bt[\cite{Stahl}]\label{StahlDomThm}
Let $f$ be a function satisfying the conditions of Note \ref{PadAssumptionsNote}  then there exists a domain $\Omega\subset \hat{\CC}$ unique up to capacity such that $\infty\in \Omega$ and 
\begin{enumerate}
    \item Any near-diagonal Padé sequence $\{[n_j/m_j]\}_{j\in \NN}$ for $f$ will converge in capacity on $\Omega$.
    \item If $\widetilde{\Omega} \supset \Omega $ and $cap(\widetilde{\Omega}\setminus \Omega)>0$ then no near-diagonal Padé sequence converges in capacity to $f$ on $\widetilde{\Omega}$.
\end{enumerate}
\et
\begin{Note}\label{SimplyConDom}
If the domain $\Omega$ from Theorem \ref{StahlDomThm} is simply connected then taking the biholomorphic map guaranteed by the Riemann mapping theorem  $\psi:\Omega \to \DD$ such that $\psi(\infty)=0$ and $\psi'(\infty)>0$ then the function $G_\Om$ from {\rm{\cite{Stahl}}} is exactly $\log |\psi|$.  So in some sense, the Padé approximants generate the conformal map as well as the maximal domain of single-valuedness.  Additionally, the next theorem allows one to extract an approximation of the conformal map from the function $G_\Om$. 
\end{Note}

Reminiscent of measure theory, we adopt the term quasi everywhere (qu.e. )to mean a statement about a set is true except possibly on subsets of capacity zero.  
\begin{Definition}[\cite{Stahl}]
A sequence of functions $\{f_n\}_{n\in\NN}$ defined on a set $\Omega $ is said to converge in capacity on $\Omega$ to a function $f$ if for every $\varepsilon>0$ and every $K\subset\subset \Omega$

\begin{equation}
    \lim_{n\to \infty} cap\{z\in K:|f_n(z)-f(z)|> \varepsilon\}=0
\end{equation}

\end{Definition}

\begin{Remark}
This set function as defined above is not generally a measure however, it is monotonic with respect to set inclusion {\rm{\cite{SaffLog}}}.  
\end{Remark}

From here we denote the unbounded component of $\CC\setminus K$ as $\Omega$ and define the Green's function with pole at $\infty$ on $\Omega$.

\begin{equation}
    g_\Omega(z,\infty) \coloneqq  V_K-U^{\mu_K}(z)
\end{equation}

There are three properties that uniquely define the above function for a given $\Omega$
\begin{enumerate}\label{GreenFuncPropert}
    \item $g_\Omega(z,\infty) $ is harmonic in $\Omega \setminus\{\infty\}$ and is bounded above and below away from every neighborhood of $\infty$. 
    \item $g_\Omega(z,\infty) -\log{|z|}=O(1)$ for $z\to \infty$
    \item $g_\Omega(z,\infty) \to 0$ as $z\to \zeta$ for qu.e. $\zeta \in \partial \Omega$. 
\end{enumerate}
Finally, we define  
\begin{equation}\label{CapG}
    G_\Omega(z)={\rm{e}}^{-g_\Omega(z,\infty)}
\end{equation}
This function carries information about the rate of convergence in capacity of the near-diagonal Padé sequence to the given function.    We can verify that $0\leq G_\Omega(z)<1$ for $z\in \Omega$ and $G_\Omega(z)>0$ for $z\in \Omega \setminus\{\infty\}$ \cite{Stahl}.

\bt [\cite{Stahl}]\label{StahlMainThm}
Let $f$ be a function satisfying the conditions of Note \ref{SimplyConDom} and $\Om$ be its maximal domain of single-valuedness from Theorem \ref{StahlDomThm}.  Furthermore, let $G_\Om$ be the corresponding Greens function \eqref{CapG} defined by \eqref{CapG} and $\NDP$ a near-diagonal Padé sequence.  Then 
\begin{enumerate}
    \item for every $\varepsilon>0$ and $V \subset \subset \Om\setminus \{\infty\}$
    $$\lim_{j\to \infty} cap\{z\in V: |f(z)-[n_j,m_j](z)|>(G_\Om(z)+\varepsilon)^{n_j+m_j}\}=0$$
    \item if $f$ has branch points, then for every $0<\varepsilon\leq \inf_{z\in V} G_\Om(z)$ and $V \subset \subset \Om\setminus \{\infty\}$ we have 
    $$\lim_{j\to \infty} cap\{z\in V: |f(z)-[n_j,m_j](z)|<(G_\Om(z)-\varepsilon)^{n_j+m_j}\}=0$$
\end{enumerate}
\et

Both parts of Theorem \ref{StahlMainThm} show that the near-diagonal Padé sequences converge in capacity in the domain established by Theorem \ref{StahlDomThm} and moreover, this convergence is geometric.  However, the major distinction between the branched and non-branched cases is the existence of a continuum in the compliment of $\Omega$.  This continuum which is composed of a union of open analytic arcs can be shown to exist if and only if the given function has branch points \cite{Stahl}.

\subsection{Dyadic Expansion of $\mathrm{Ei}$}\label{DyEiSec}
Using the general techniques developed in \cite{GRA} where we obtained \eqref{eq:EidsumP1} a geometrically convergent rational approximation of   ${\rm Ei}$ that is valid in the cut plane $\CC\setminus -i\overline{\RR^+}$.  For a review of this special function see \S\ref{EidetailsRatApp} or for the full theory see \cite{GRA}.

\begin{equation}\label{eq:EidsumP1}
   \mathrm{e}^{-x}{\rm Ei}^+(x)= -\sum_{m=1}^{\infty}\frac{\Gamma(m)}{2^m}\frac{1}{(y)_m}+\sum_{k=1}^\infty\sum_{m=1}^\infty\frac{\Gamma(m) \mathrm{e}^{- i\pi/2^{k}}}{(1+\mathrm{e}^{- i\pi/2^{k}})^m}\frac{1}{(2^ky)_m} \ \ \ \ \ \ (y=-i x/\pi)\qquad 
\end{equation}
where 
  \begin{equation}\label{PochRat}
 (x)_k:=x(x+1)\cdots (x+k-1)=\frac{\Gamma(x+k)}{\Gamma(x)}
  \end{equation}
 is known as the Pochhammer symbol, or rising factorial. 

For approximations we need truncated series and estimates of the remainder. 
We use the following notation 
 \begin{equation}
   \label{rhon}
   \rho_{n+1,0}
   (z,x):=(1-z)\sum_{k= n+1}^\infty z^k\frac{k!}{(x)_{k+1}}
   \end{equation}

   and
   
For $x\notin (-\infty,0]\beta$ we have
\begin{equation}
\label{RNEstThm}
    R_N(t,x; \beta
) =\frac{1}{\beta}\int_{\Gamma_{c}} \mathrm{e}^{-xq/\beta} \rho_{N}(q/\beta,1+t {\rm{e}}^{i\theta};\beta)dq\\\ 
\end{equation}
with $\rho_{N}$ 
given by \eqref{eq:remnbeta} and $\beta \in \CC\setminus \{0\}$ a free complex parameter that allows us to place the ray of poles 
in the dyadic expansion along a direction of our choice c.f. \cite{GRA}.  In \eqref{eq:EidsumP1} we have taken $\beta=\pi i$ and one checks that 
the poles all lie along the ray $(-\infty,0]i$.  The contour  \( \Gamma_{c} \) is chosen so that along it we have: (a)  $\Re(xq/\beta)>0$ for large $q$ and (b) the function $\rho_N$ (defined by \eqref{eq:remnbeta} with $1+t {\rm{e}}^{i\theta}=s$) is analytic. More precisely,
given \( x\in \mathbb{C} \setminus \beta(-\infty,0] \) and any  
\( c\in \left(0,\min\left\{1,\frac{1}{2}\operatorname{dist}(x,\beta(-\infty,0])\right\}\right) \),  define

\begin{equation}
\label{RemContGam}
    \Gamma_{c}=\begin{cases}
        \RR^+ & \Re(x/\beta)>c\\[6pt]
       \mathrm{e}^{-i\pi/2}\RR^+ &   \Re(x\mathrm{e}^{-i\pi/2}/\beta)>c\\[6pt]
        \mathrm{e}^{i\pi/2}\RR^+ &   \Re(x\mathrm{e}^{i\pi/2}/\beta)>c
        
    \end{cases}
\end{equation}  
For large $N$, $x\in\Omega_c$ and $c>0$,
\begin{equation}
\label{eq:singrem}
    |R_N(t, x;\beta
    )|\le
     \begin{cases}
      \frac1{c_02^{N-1}\Re (x/\beta)}& \Re(x/\beta)>c\\[10pt]
    \frac{1}{c_02^{N-1}\Re(x\mathrm{e}^{-i\pi/2}/\beta)}& x\in\{\Re(x\mathrm{e}^{-i\pi/2}/\beta)> c\}\cap\{\Re (x/\beta)\leq c\}\\[10pt]
    \frac{1}{c_02^{N-1}\Re(x\mathrm{e}^{i\pi/2}/\beta)}& x\in\{\Re(x\mathrm{e}^{i\pi/2}/\beta)> c\}\cap\{\Re (x/\beta)\leq c\}
    \end{cases}
 \end{equation} 

 \begin{equation}
  \label{eq:remnbeta}
    \rho_{n+1}(p,s;\beta
    ):=\frac{\beta}{2^n}\left( \frac 1{\beta(s-p)/2^n} +\frac{  \mathrm{e}^{-{\beta s}/{2^n}}}{\mathrm{e}^{-\beta {s}/{2^n}}-\mathrm{e}^{-\beta{p}/{2^n}}}  \right).
    \end{equation}

Writing the series \eqref{eq:EidsumP1} as a sum with remainder we have
\begin{equation}
\label{dyadicEi}
 \mathrm{e}^{-x}{\rm Ei}^+(x)=-\sum_{m=1}^n \frac{\Gamma(m)}{2^m\left(\frac x{i\pi}\right)_m} -\rho_{n,0}(x)+\sum_{k=1}^{N-1} \left[ \sum_{m=1}^\ell \frac{\Gamma(m)\,\mathrm{e}^{- i\pi/2^{k}}}{\left(1+\mathrm{e}^{- i\pi/2^{k}}\right)^{m} \,  \left(\frac {2^kx}{i\pi}\right)_m}+\rho_{\ell,k}(x)\right]+R_N(0,x;\pi i)
 \end{equation}
 where the remainders are 
\begin{equation}
\label{remEiRatApp}
\rho_{n,0}(x)=\rho_{n,0}
\left(\frac12,\frac x{i\pi}\right),\ \ \ \ \ \rho_{\ell,k}(x)= \rho_{\ell,0
}\left(\frac{e_k}{1+e_k},\frac {2^kx}{i\pi}\right),\ \ \ \ \text{ where }e_k={\rm{e}}^{i\pi/2^k} 
\end{equation}
which are defined by \eqref{rhon}

\begin{Proposition}\label{P2} 
 
   (i) For fixed $x\in \CC\setminus -i\overline{\RR^+}$ and large $n$,  $\rho_{n,0}(x)=O(2^{-n}n^{-\Im x/\pi})$. For fixed $n$ and large $x$, $\rho_{n,0}(x)=O(x^{-n})$. 

(ii) For large $l$, 
fixed $k$ and $x\in \Omega_c$, $\rho_{\ell,k}(x)=O(  |1+e_k^{-1}|^{-\ell}\ell^{-2^k\Im x/\pi})$. For fixed $\ell$ and large $2^kx$, $\rho_{\ell,k}(x)=O((2^{k}x)^{-\ell})$.

(iii) For large $N$ and fixed $x$ such that $\Re(x/ \pi i)>c$ we have 
  $|R_N(0,x;\pi i)|\leq \frac{1}{c_0 2^{N-1} \Re(x/\pi i)}$.  For the other regions $\Re(x{\rm{e}}^{-i \pi}/\pi )> c$ and $\Re(x/\pi )> c$ similar estimates follow from \eqref{eq:singrem}.
\end{Proposition}

   \subsection{More about Ei}\label{EidetailsRatApp}
  
We observe that ${\rm{E_1}}$ can be expressed in a manner conducive to analytic continuation across the discontinuity on $\mathbb{R}_-$: through elementary substitutions, we can rewrite it as follows, for $z>0$,
  $${\rm{E_1}}(z):=\int_z^\infty\frac{{\rm{e}}^{-t}}{t}\, dt= {\rm{e}}^{-z}\int_0^\infty \frac{{\rm{e}}^{-pz}}{1+p}\, dp$$
 
 It is interesting to note the relation  6.2.4 \cite{nist}
 $${\rm{E_1}}(z)={\rm Ein}(z)-\ln z-\gamma,\ \ \ \ \ \ \text{where }\ \ \ {\rm Ein}(z)=\int_0^z\frac{1-{\rm{e}}^{-t}}{t}\, dt$$

Given that ${\rm{Ein}}$ is an entire function, and $\ln$ is defined with the standard branch for $z > 0$, it follows that during analytic continuation across $\mathbb{R}_-$, the function ${\rm{E_1}}$ acquires a $2\pi i$. Consequently, $\mathbb{R}_-$ constitutes a Stokes line.
 
 For us it is convenient to place this Stokes line along $\RR_+$, so we work with $-E_1(-z)$:
 $$-{\rm{E_1}}(-z)= -{\rm{e}}^{z}\int_0^\infty \frac{{\rm{e}}^{pz}}{1+p}\, dp,\ \ (z<0)$$
 which analytically continued to the first quadrant yields our Ei$^+$ defined by \eqref{LiEiRatAp}.

 \begin{equation}
\label{LiEiRatAp}
{\rm Ei}^+(x)=\mathrm{e}^{x}\int_0^{\infty {\rm{e}}^{i0-}}\frac{\mathrm{e}^{-px}}{1-p}dp
\end{equation}
 
 Note the structure of the branch point at $0$:
 
 $${\rm Ei}^+(z)=-{\rm Ein}(-z)+\ln (-z)+\gamma$$
 where $\ln (-z)$ has the usual brach for $z<0$ and then it is analytically continued on the Riemann surface of the log (and ${\rm Ein}(-z)$ is entire).
 
 For ${\rm{e}}^{-x}$Ei$^+(x)$, $\RR^+$ is a Stokes ray and the two sides of $i\RR^-$ are antistokes lines. The behavior of ${\rm{e}}^{-x}$Ei$^+(x)$ is oscillatory when $\arg(x)=-\pi/2$ and it is given by an asymptotic series when $\arg(x)=3\pi/2$.

\section{\textbf{Acknowledgements}}  I would like to thank Dr. Ovidiu Costin for his patience and dedication to his students.  Without him, this project would not have been possible.

\end{document}